\documentstyle[twoside,amstex]{article}
\input amssym.def
\input amssym.tex
\def\bc{\begin{center}}
\def\ec{\end{center}}
\def\no{\noindent}
\def\hang{\hangindent\parindent}
\def\textindent#1{\indent\llap{[#1]\enspace}\ignorespaces}
\def\re{\par\hang\textindent}
\begin{document}
\thispagestyle{empty} \vspace*{3 true cm} \pagestyle{myheadings}
\markboth {\hfill {\sl M.S. Abdolyousefi, R. Bahmani and H. Chen}\hfill}
{\hfill{\sl Elementary Matrix Reduction over Locally Stable Rings}\hfill} \vspace*{-1.5 true cm} \bc{\large\bf Elementary Matrix Reduction over Locally Stable Rings}\ec

\vskip6mm
\bc{{\bf Marjan
Sheibani Abdolyousefi}\\[3mm]
Faculty of Mathematics, Statistics and Computer Science\\
Semnan University, Semnan, Iran\\
m.sheibani1@@gmail.com}\ec

\bc{{\bf Rahman Bahmani}\\[3mm]
Faculty of Mathematics, Statistics and Computer Science\\
Semnan University, Semnan, Iran\\
rbahmani@@semnan.ac.ir}\ec

\bc{{\bf Huanyin Chen $^*$}\\[3mm]
Department of Mathematics, Hangzhou Normal University\\
Hangzhou 310036, China\\
huanyinchen@@aliyun.com}\ec

\begin{figure}[b]
\vspace{-3mm}
\rule[-2.5truemm]{5cm}{0.1truemm}\\[3mm]
{\footnotesize 2010 Mathematics Subject Classification. 13F99, 13E15, 06F20.
$^*$ Corresponding author.}
\end{figure}

\vskip10mm
\begin{minipage}{120mm}
\no {\bf Abstract:} A commutative ring is locally stable (has neat range 1) provided that for any $a,b\in R$ such that $aR+bR=R$, there exists a
$y\in R$ such that $R/(a+by)R$ has stable range 1 (is clean). For a B$\acute{e}$zout ring $R$, we prove that $R$ is an elementary divisor ring if and only
if $R$ is a locally stable ring if and only if $R$ has neat range 1. This completely determines when a B$\acute{e}$zout ring is an elementary divisor ring. In particular, we answer two open problems of Zabavsky and Pihura. Lastly, we prove that every locally stable ring is strongly completable.
\end{minipage}

\vskip15mm \bc{\bf 1. Introduction}\ec

\vskip4mm \no Throughout this paper, all rings are commutative
with an identity. A matrix $A$ (not necessarily square) over a ring
$R$ admits diagonal reduction if there exist invertible matrices
$P$ and $Q$ such that $PAQ$ is a diagonal matrix $(d_{ij})$, for
which $d_{ii}$ is a divisor of $d_{(i+1)(i+1)}$ for each $i$. A
ring $R$ is called an elementary divisor ring provided that every
matrix over $R$ admits a diagonal reduction. A ring is a
B$\acute{e}$zout ring if every finitely generated ideal is
principal. Obviously, $\{ ~\mbox{elementary divisor rings}~\}\subsetneq \{ ~~\mbox{B$\acute{e}$zout rings}~\}$ (cf. [11, Theorem 1.2.11 and Corollary 2.1.3]). Gillman and Heriksen constructed an example of a B$\acute{e}$zout ring which is not an elementary divisor ring. An attractive problem is to investigate various conditions under which a B$\acute{e}$zout ring is an elementary divisor ring (cf. [8], [10], [12] and [14]).

A ring $R$ has stable range 1 if for any $a,b\in R$ such that $aR+bR=R$, there exists a $y\in R$ such that $a+by\in U(R)$. For instance, every clean ring (i.e., the ring in which every element is the sum of an idempotent and a unit) has stable range 1. For general theory of such rings, we refer the reader to the book [3]. The main purpose of this paper is to characterize elementary divisor rings by means of a type of stable-like conditions. We call a ring $R$ is locally stable provided that for any $a,b\in R$ such that $aR+bR=R$, there exists a $y\in R$ such that $R/(a+by)R$ has stable range 1.

In Section 2, we investigate elementary properties of locally stable rings, which will be used repeatedly in the sequel.
We shall see that locally stable rings contain almost all known rings of interest. In Section 3, we explore when a B$\acute{e}$zout ring is an
elementary divisor ring. Surprisingly, we prove that a ring is an elementary divisor ring if and only if it is a locally stable B$\acute{e}$zout ring if
and only if $R$ is a B$\acute{e}$zout ring of neat range 1. Here, a ring $R$ has neat range 1 provided that for any $a,b\in R$ such that $aR+bR=R$,
there exists a $y\in R$ such that $R/(a+by)R$ is clean [12]. This completely determines when a B$\acute{e}$zout ring is an elementary ring. In particular, we answer two open problems of Zabavsky and Pihura [15, Open Problems].
Finally, we prove that every locally stable ring is strongly completable. Many known results are thereby
generalized to much wider class of rings, e.g. [4, Theorem 14], [8, Theorem 3.7], [9, Theorem ], [11, Theorem 1.2.13 and Theorem 1.2.21] and [12, Theorem 32].

We shall use $J(R)$ and $U(R)$ to denote the Jacobson radical of $R$ and the set of all units in $R$, respectively. $M_n(R)$ denotes
the ring of all $n\times n$ matrices over $R$, and $GL_n(R)$ stands for the $n$-dimensional general linear group of $R$.

\vskip15mm\bc{\bf 2. Locally Stable Rings}\ec

\vskip4mm We call an element $a\in R$ is stable if $R/aR$ has stable range 1. Thus, a ring $R$ is locally stable if and only if for any $a,b\in R$ such that $aR+bR=R$ there exists a $y\in R$ such that $a+by\in R$ is stable. Clearly, every ring of stable range 1 is locally stable. But the converse is not true, e.g., ${\Bbb Z}$. The purpose of this section is to investigate elementary properties and provide various locally stable rings. We begin with

\vskip4mm \hspace{-1.8em} {\bf Proposition 2.1.}\ \ {\it Let $R$ be a ring. Then the following are equivalent:}
\vspace{-.5mm}
\begin{enumerate}
\item [(1)]{\it $R$ is locally stable.}
\vspace{-.5mm}
\item [(2)]{\it For any $c\in R$, $\overline{u}\in U(R/cR)\Longrightarrow ~u-a\in cR$ for some stable element $a\in R$.}
\end{enumerate}
\vspace{-.5mm} {\it Proof.}\ \  $(1)\Rightarrow (2)$ Let $\overline{u}\in U(R/cR)$. Then there exists some $\overline{v}\in R/cR$ such that $\overline{u}\overline{v}=\overline{1}$. So we can find some $x\in R$ such that $uv+cx=1$. Then $uR+cR=R$. As $R$ is locally stable,
$R/(u+cy)R$ has stable range 1 for a $y\in R$. Let $a=u+cy$, then $a\in R$ is stable and $u-a\in cR$.

$(2)\Rightarrow (1)$ Let $aR+bR=R$ for some $a,b\in R$. Then there exist $x,y\in R$ such that $ax+by=1$. We have $ax+bR=1+bR$, that implies
$\overline{a}\in U(R/bR)$. By $(2)$, there exists some stable element $y\in R$ such that $a-y\in bR$. Hence $y=a+br$ for some $r\in R$ and $R/(a+br)R$
has stable range 1. This
shows that $R$ is locally stable.\hfill$\Box$

\vskip4mm Following McGovern, a ring $R$ has almost stable range 1 if every proper homomorphism of $R$ has stable range 1. For instance, rings having stable range 1 (e.g., clean rings, local rings, semilocal rings, etc), adequate rings, neat rings, semi-clean rings all have almost stable range 1 (cf. [8]).

\vskip4mm \hspace{-1.8em} {\bf Corollary 2.2.}\ \ {\it If $R$ has almost stable range 1, then $R$ is locally stable.}
\vskip2mm\hspace{-1.8em} {\it Proof.}\ \ Let $\overline{u}\in U(R/cR)$ for any $c\in R$. Then $u\neq cr$ for any $r\in R$. Let $a=c-u$,
then $R/aR$ has stable range 1, since $a\neq 0$ and $R$ has almost stable range 1. Hence by Proposition 2.1. $R$ is locally stable.\hfill$\Box$

\vskip4mm But the converse is not true as the following shows.

\vskip4mm \hspace{-1.8em} {\bf Example 2.3.}\ \ {\it Let $R=\{ a_0+a_1x+a_2x^2+\cdots ~|~a_0\in {\Bbb Z}, a_1,a_2,\cdots \in {\Bbb Q}\}$. Then $R$ is locally stable, but $R$ has no almost stable range 1.}
\vskip2mm\hspace{-1.8em} {\it Proof.}\ \  As $R/{\Bbb Q}[[X]]\cong {\Bbb Z}$ is a homomorphic image of $R$ with stable range 2, then $R$ does not have
almost stable range 1. Now let $fR+gR=R$ for some $f,g\in R$, then $f(0){\Bbb Z}+
g(0){\Bbb Z}={\Bbb Z}$. Clearly, we can find a $y\in {\Bbb Z}$ such that $f(0)+g(0)y\neq 0$, we see that ${\Bbb Z}/(f(0)+g(0)y){\Bbb Z}$ is finite, and so it has stable range 1.
We shall show that $R/(f+gy)R$ has stable range 1. Given $\overline{cR+dR}=R/(f+gy)R$, then $\overline{c(0){\Bbb Z}+d(0){\Bbb Z}}={\Bbb Z}/(f(0)+g(0)y){\Bbb Z}$. Thus,
we can find a $z\in {\Bbb Z}$ such that $\overline{c(0)+d(0)z}\in U\big({\Bbb Z}/(f(0)+g(0)y){\Bbb Z}\big)$.
Write $\overline{(c(0)+d(0)z)h}=\overline{1}$, and so $1-(c(0)+d(0)z)h\in (f(0)+g(0)y){\Bbb Z}$. We check that
$$1-(c(0)+d(0)z)h\in (f(x)+g(x)y)p+x{\Bbb Q}[[x]].$$ Thus, we have a $q(x)\in {\Bbb Q}[[x]]$ such that
$$1-xq(x)-(c(0)+d(0)z)h=(f(x)+g(x)y)p.$$ Further, we can find a $q'(x)\in {\Bbb Q}[[x]]$ such that
$$1-xq(x)-xq'(x)-(c(x)+d(x)z)h=(f(x)+g(x)y)p.$$ As $1-xq(x)-xq'(x)\in U({\Bbb Q}[[x]])$, we get
$$1-(c(x)+d(x)z)h\big(1-xq(x)-xq'(x)\big)^{-1}=(f(x)+g(x)y)p\big(1-xq(x)-xq'(x)\big)^{-1}.$$ Therefore $\overline{c(x)+d(x)z}\in U(R/(f+gy)R),$ as required.
\hfill$\Box$

\vskip4mm An element $c\in R$ is adequate provided that for any $a\in R$ there exist some $r,s\in R$ such that $(1)$ $c=rs$; $(2)$ $rR+aR=R$; $(3)$ $s'R+cR\neq R$ for each non-invertible divisor $s'$ of $s$. A
B$\acute{e}$zout ring in which every nonzero element is adequate
is called an adequate ring. Following Domsha and Vasiunyk,, a ring $R$ has adequate range 1 provided that $aR+bR=R$ with $a,b\in R\Longrightarrow ~\exists y\in R$ such that $a+by\in R$ is adequate [4]. For instance, rings having stable range 1, adequate rings and VNL rings ( i.e. for any element $a$, either $a$ or $1-a$ is regular) have adequate range 1.
Following Zabavsky, a ring $R$ has neat range 1 if $aR+bR=R$ with $a,b\in R\Longrightarrow ~\exists~y\in R$ such that $R/(a+by)R$ is clean ([12]). We have

\vskip4mm \hspace{-1.8em} {\bf Example 2.4.}
\begin{enumerate}
\item [(1)] Every ring of adequate range 1 is locally stable.
\item [(2)] Every ring of neat range 1 is locally stable.
\end{enumerate}
\vspace{-.5mm} {\it Proof.}\ \ (1) Let $R$ have adequate range 1. For any $a,b\in R$ such that $aR+bR=R$, there exists a $y\in R$ such that $a+by\in R$ is adequate. In view of [13, Theorem 8], $R/(a+by)R$ has stable range 1. Thus, $R$ is locally stable, as desired.

(2) Let $R$ have neat range 1, and let $aR+bR=R$ with $a,b\in R$. Then there exists a $y\in R$ such that $R/(a+by)R$ is clean. In view of [3, Corollary 1.3.15], $R/(a+by)R$ has stable range 1, and thus $R$ is locally stable. This completes the proof.\hfill$\Box$

\vskip4mm Every homomorphic image of a locally stable range is locally stable. Further, we have

\vskip4mm \hspace{-1.8em} {\bf Proposition 2.5.}\ \ {\it Let $I$ be an ideal of $R$, and let $I\subseteq J(R)$. Then $R$ is locally stable if and only if so is $R/I$.}
\vskip2mm\hspace{-1.8em} {\it Proof.}\ \ Let $R$ be a locally stable ring and let  $I\subseteq J(R)$ be an ideal of $R$, such that
$\overline{a}\overline{R}+\overline{b}\overline{R}=\overline{R}$, where $\overline{R}=R/I$. There exist some $\overline{x},
\overline{y}\in R/I$ and $\overline{a}\overline{x}+\overline{b}\overline{y}=\overline{1}$. This shows that  $ax+by=1+r$ for
some $r\in I$. As $r\in I\subseteq J(R)$, $ax+by\in U(R)$. Hence, $aR+bR=R$. Now $R/(a+bz)R$ has stable range 1 for some $z\in R$ and we deduce
that $\overline{R}/(\overline{a}+\overline{bz})\overline{R}$ has stable range 1. Accordingly, $R/I$ has stable range 1.

Conversely assume that $R/I$ is locally stable. Let $\overline{u}\in U(R/cR)$ for some $c\in R$, then
$\overline{u+I}\in U(R/I/(c+I)R/I)$, by Proposition 2.1. There exists a stable element $a+I\in R/I$ such that
$(u+I)-(a+I)\in (c+I)R$ . So we have $u-(a+w)\in cR$ for some $w\in I$. As $R/I/(a+I)R/I$ has stable range 1 and $I\subseteq J(R)$, one easily checks that
$R/(a+w)R$ has stable range 1. Therefore $R$ is locally stable.\hfill$\Box$

\vskip4mm As an immediate consequence of Proposition 2.5, we see that $R$ is locally stable if and only if so is $R/J(R)$. Moreover, we now derive

\vskip4mm \hspace{-1.8em} {\bf Corollary 2.6.}\ \ {\it A ring $R$ is locally stable if and only if so is $R[[X]]$.}
\vskip2mm\hspace{-1.8em} {\it Proof.}\ \ Let $\varphi:R[[X]]\rightarrow R$ be a morphism, that is defined by $\varphi(f(X))=f(0)$ for any $f(X)\in R[[X]]$.
It is obvious that $ker(\varphi
)\subseteq J(R[[X]])$ and $R[[X]]/ker(\varphi)\cong R$. Now by the Proposition 2.5. $R[[X]]$ is locally stable if and only if $R$ is locally stable.    \hfill$\Box$

\vskip4mm \hspace{-1.8em} {\bf Proposition 2.7.}\ \ {\it Let $\{ R_i\}$ be a family of rings. Then the direct product $\prod\limits_{i}R_i$ of rings is locally stable if and only if each $R_i$ is locally stable.}
\vskip2mm\hspace{-1.8em} {\it Proof.}\ \ $\Longrightarrow $ Suppose that $T:=\prod\limits_{i}R_i$ is locally stable.
Given $aR_1+bR_1=R_1$ for some $a,b\in R_1$, we have $(a,1,1,\cdots )T+(b,0,0,\cdots )T=T$. By hypothesis, there exist $(y,s_2,s_3,\cdots )\in T$ such that $T/((a,1,1,\cdots )+(b,0,0,\cdots )(y,s_2,s_3,\cdots ))T$ has stable range 1, but we have $T/((a,1,1,\cdots )+(b,0,0,\cdots )(y,s_2,s_3,\cdots ))T\cong R_1/(a+by)R_1$. Therefore $R_1/(a+by)R_1$ has stable range 1, and then $R_1$ is locally stable. We can see in the same manner that $R_i (i\neq 1)$ is locally stable.

$\Longleftarrow $ Suppose that each $R_i$ is locally stable and $(a_1,a_2,\cdots )T+(b_1,b_2,\cdots )T=T$, where $T=\prod\limits_{i}R_i$. Then we have $a_iR_i+b_iR_i=R_i$. Hence, there exist $y_i\in R_i$ such that $R_i/(a_i+b_iy_i)R$ has stable range 1. One easily checks that
$$T/((a_1,a_2,\cdots )+(b_1,b_2,\cdots )(y_1,y_2,\cdots ))T\cong \prod\limits_{i}R_i/(a_i+b_iy_i)R_i,$$ and therefore $T/((a_1,a_2,\cdots )+(b_1,b_2,\cdots )(y_1,y_2,\cdots ))T$ has stable range 1. This completes the proof.\hfill$\Box$

\vskip4mm \hspace{-1.8em} {\bf Corollary 2.8.}\ \ {\it Let $e\in R$ be an idempotent of $R$. Then $R$ is locally stable if and only if $eRe$ and $(1-e)R(1-e)$ are locally stable.}
\vskip2mm\hspace{-1.8em} {\it Proof.}\ \ For any idempotent $e\in R$ we have, $R\cong eRe\oplus (1-e)R(1-e)$, as $R$ is commutative. Now the result follows from the Proposition 2.7.
\hfill$\Box$

\vskip4mm \hspace{-1.8em} {\bf Theorem 2.9.}\ \ {\it Let $R$ be a ring. If $a$ or $1-a$ is stable for all $a\in R$, then $R$ is locally stable.}
\vskip2mm\hspace{-1.8em} {\it Proof.}\ \ Given $aR+bR=R$ with $a,b\in R$, then $ax+by=1$ for some $x,y\in R$. Hence, $a(x-y)+(a+b)y=1$. By hypothesis, $R/a(x+y)R$ or $R/(a+b)yR$ has stable range 1.

Case 1. $R/a(x+y)R$ has stable range 1. Clearly, $R/aR\cong R/(a(x+y))R/aR/(a(x+y))R$. It follows that $R/aR$ has stable range 1. That is, $R/(a+b\cdot 0)R$ has stable range 1.

Case 2. $R/(a+b)yR$ has stable range 1. As in the proof in Case I, we see that $R/(a+b\cdot 1)R$ has stable range 1.

Accordingly, we have a $z\in R$ such hat $R/(a+bz)R$ has stable range 1. Therefore, $R$ is locally stable.\hfill$\Box$

\vskip4mm Following Domsha and Vasiunyk [4], a ring $R$ is a locally adequate ring if for every $a\in R$, $a$ or $1-a$ is adequate. For instance, NJ-rings, i.e., rings in which every element not in the Jacobson radicals is regular.

\vskip4mm \hspace{-1.8em} {\bf Corollary 2.10.}\ \ {\it Every locally adequate ring is locally stable.}
\vskip2mm\hspace{-1.8em} {\it Proof.}\ \ Let $a\in R$. Then $a$ or $1-a$ is adequate. In view of [13, Theorem 8], $R/aR$ or $R(1-a)R$ has stable range 1. That is, $a$ or $1-a$ is stable. Therefore the proof is true, by Theorem 2.9.\hfill$\Box$

\vskip15mm\bc{\bf 3. Elementary Divisor Rings}\ec

\vskip4mm Recall that a ring $R$ has stable range 2 provided that $aR+bR+cR=R$ with $a,b,c\in R\Longrightarrow ~\exists y,z\in R$ such that $(a+cy)R+(b+cz)R=R$.

\vskip4mm \hspace{-1.8em} {\bf Lemma 3.1.}\ \ {\it Every locally stable ring has stable range 2.}
\vskip2mm\hspace{-1.8em} {\it Proof.}\ \ Let $R$ be a locally stable ring. Suppose that $aR+bR+cR=R$ with $a,b,c\in R$.
Then there exist $y,z\in R$ such that $w:=a+by+cz$ and $R/wR$ has stable range 1. Write $ak+bp+cq=1$ with $k,p,q\in R$. Then
$(a+by+cz)k+b(p-yk)+c(q-zk)=1$. Hence, $wR+bR+cR=1$, and so
$\overline{b}(R/wR)+\overline{c}(R/wR)=R/wR$. Thus, we can find some $d\in R$ such that
$\overline{b+cd}\in U(R/wR)$; whence, $wR+(b+cd)R=R$. Write $ws+(b+cd)t=1$. Then
$\big(a+(b+cd)y+c(z-dy)\big)s+(b+cd)t=1$, and so $\big(a+c(z-dy)\big)s+(b+cd)(t+yt)=1$.
Accordingly, $\big(a+c(z-dy)\big)R+(b+cd)R=R$, and thus yielding the result.
\hfill$\Box$

\vskip4mm \hspace{-1.8em} {\bf Lemma 3.2 ([8, Theorem 1.1] and [10, Theorem 2.5]).}\ \ {\it Let $R$ be a Hermite ring. Then the following are
equivalent:}
 \vspace{-.5mm}
\begin{enumerate}
\item [(1)] {\it $R$ is an elementary divisor ring.} \vspace{-.5mm}
\item [(2)] {\it Every matrix $\left(
\begin{array}{cc}
a&0\\
b&c
\end{array}
\right)\in M_2(R)$ with $aR+bR+cR=R$ admits an elementary reduction.}
\end{enumerate}

\vskip4mm We are now ready to prove the following.

\vskip4mm \hspace{-1.8em} {\bf Theorem 3.3.}\ \ {\it If $R$ is locally stable, then $R$ is an elementary divisor ring if and only if $R$ is a Bezout ring.}
\vskip2mm\hspace{-1.8em} {\it Proof.}\ \ $\Longrightarrow$ This is obvious.

$\Longleftarrow$ By virtue of Lemma 3.1, $R$ has stable range 2. It follows from [11, Theorem 2.1.2] that $R$ is a Hermite ring. Suppose that
$A=\left(
\begin{array}{cc}
a&0\\
b&c
\end{array}
\right)\in M_2(R)$ with $aR+bR+cR=R$. In terms of Lemma 3.2, it suffices to check that $A$ admits an elementary reduction.

Write $ax+by+cz=1$. Then $bR+(ax+cz)R=R$. As $R$ is locally stable, there exists some $t\in R$ such that $v=b+(ax+cz)t$ and $R/vR$ has stable range 1. We have $$\left(
\begin{array}{cc}
1&0\\
xt&1
\end{array}
\right) \left(
\begin{array}{cc}
a&0\\
b&c
\end{array}
\right)\ \left(
\begin{array}{cc}
1&0\\
zt&1
\end{array}
\right)=\left(
\begin{array}{cc}
a&0\\
v&c
\end{array}
\right).$$ Since $R$ is a Hermite ring, there exists $P\in GL_2(R)$ such that $(v,c)P=(0,c^{\prime})$, and so $\left(
\begin{array}{cc}
a&0\\
v&c
\end{array}
\right)P=\left(
\begin{array}{cc}
a^{\prime}&b^{\prime}\\
0&c^{\prime}
\end{array}
\right)$. It is easily seen that $vR\subseteq c^{\prime}R$ and $a^{\prime}R+b^{\prime}R+c^{\prime}R=R$. Since $R/vR$ has stable range 1, and so $R/c^{\prime}R$ has stable range 1, as
$R/c'R\cong (R/vR)/(c'R/vR)$. Also $\overline{a^{\prime}}(R/c^{\prime}R) + \overline{b^{\prime}}(R/c^{\prime}R)=R/c^{\prime}R$, hence there exists some $w\in R$ such that $\overline{b^{\prime}+a^{\prime}w}\in U(R/c^{\prime}R)$, and so $(\overline{b^{\prime}+a^{\prime}w)p}=\overline{1}$, for some $p\in R$, and then $(b'+aw)p+c'q=1$. One easily checks that
$$\begin{array}{ll}
&\left(
\begin{array}{cc}
&1\\
1&
\end{array}
\right)\left(
\begin{array}{cc}
c'&-(b'+aw)\\
p&q
\end{array}
\right)\left(
\begin{array}{cc}
a^{\prime}&b'\\
&c^{\prime}
\end{array}
\right)\left(
\begin{array}{cc}
1&w\\
&1
\end{array}
\right)\left(
\begin{array}{cc}
1&\\
-pa'&1
\end{array}
\right)\left(
\begin{array}{cc}
&1\\
1&
\end{array}
\right)\\
=&\left(
\begin{array}{cc}
1&\\
&a'c'
\end{array}
\right).
\end{array}$$
As $det\left(
\begin{array}{cc}
c'&-(b'+aw)\\
p&q
\end{array}
\right)=1$, we see that $\left(
\begin{array}{cc}
c'&-(b'+aw)\\
p&q
\end{array}
\right)\in GL_2(R)$. Thus, we prove that
$\left(
\begin{array}{cc}
a^{\prime}&b^{\prime}\\
&c^{\prime}
\end{array}
\right)$ admits a diagonal reduction. Therefore $\left(
\begin{array}{cc}
a&0\\
b&c
\end{array}
\right)$ admits a diagonal reduction, and so $R$ is an elementary divisor ring.\hfill$\Box$

\vskip4mm \hspace{-1.8em} {\bf Corollary 3.4 ([8, Theorem 3.7]).}\ \ If $R$ has almost stable range 1, then $R$ is an elementary divisor ring if and only if $R$ is a B$\acute{e}$zout ring.

\vskip4mm The following corollary extends [4, Theorem 14] to rings which may have many zero-divisors.

\vskip4mm \hspace{-1.8em} {\bf Corollary 3.5.}\ \ If $R$ has adequate range 1, then $R$ is an elementary divisor ring if and only if $R$ is a B$\acute{e}$zout ring.

\vskip4mm \hspace{-1.8em} {\bf Lemma 3.6.}\ \ {\it Let $R$ be a B$\acute{e}$zout ring, and let $0\neq c\in R$. Consider the following conditions.}
\vspace{-.5mm}
\begin{enumerate}
\item [(1)]{\it For any $a,b\in R$ such that $aR+bR=R$, there exist $r,s\in R$ such that $c=rs$ and $rR+sR=rR+aR=sR+bR=R$.}
\vspace{-.5mm}
\item [(2)]{\it $R/cR$ is clean.}
\end{enumerate}
Then $(1)\Rightarrow (2)$. If $R$ is a domain, then $(1)\Longleftrightarrow (2)$.
\vskip2mm\hspace{-1.8em} {\it Proof.}\ \ $(1)\Rightarrow (2)$ Given $\overline{a}(R/cR)+\overline{b}(R/cR)=R/cR$, we have $aR+bR+cR=R$. Write $bR+cR=dR$ for a $d\in R$. Then $aR+dR=R$. By hypothesis,
there exist $r,s\in R$ such that $c=rs$ and $rR+sR=rR+aR=sR+dR=R$. Since $rR+sR=R$, there exist $u,v\in R$ such that $ru+sv=1$. Hence, $\overline{r^2u}=\overline{r}$ and $\overline{s^2v}=\overline{s}$ in $R/cR$. Let $\overline{e}=\overline{sv}$. Then $\overline{e}\in R/cR$ is an idempotent.
As $rR+aR=R$, we see that $reR+eaR=eR$, and so $\overline{e}=\overline{a\beta e}$ for some $\beta\in R$. Since
$sR+dR=R$, we see that $s(1-e)R+(1-e)dR=(1-e)R$. Thus, $(1-e)dR\equiv (1-e)R~(~mod~cR)$. It follows that $(1-e)bR\equiv (1-e)R~(~mod~cR)$, and thus $\overline{1-e}=\overline{b\gamma (1-e)}$ for some $\gamma\in R$. Therefore we have an idempotent $\overline{e}\in \overline{a}(R/cR)$ such that $\overline{1-e}\in \overline{b}(R/cR)$. This shows that $R/cR$ is clean.

Suppose that $R$ is a domain. It will suffice to prove $(2)\Rightarrow (1)$. Given $aR+bR=R$ with $a,b\in R$, we have $\overline{a}\overline{R}+\overline{b}\overline{R}=\overline{R}$, where $\overline{R}=R/cR$. Then we can find an idempotent $\overline{e}\in \overline{a}\overline{R}$ such that $\overline{1-e}\in \overline{b}\overline{R}$.  Write $e-ap=cq, 1-e-b\alpha=c\beta$ and $e(1-e)=ct$. Since $R$ is a B$\acute{e}$zout domain, there exists a $d\in R$ such that
$eR+cR=dR$. Write $e=de_0,c=dc_0$ and $d=ex+cy$. Then $d(e_0x+c_0y-1)=0$, and so $e_0x+c_0y=1$.
As $de_0(1-e)=ct=dc_0t$, we get $e_0(1-e)=c_0t$. Thus, $e_0x(1-e)+c_0y(1-e)=1-e$; hence, $e+c_0j=1$ for a $j\in R$. Take $r=c_0$ and $s=d$. Then $c=rs$. As $e=de_0=se_0$, we get $se_0+rj=1$; whence, $rR+sR=R$. Since $e=ap+cq$, we have $ap+cq+c_0j=1$, and then $ap+r(dq+j)=1$. This implies that $rR+aR=R$.
Since $c\beta+b\alpha+e=1$, we get $d(c_0\beta+e_0)+b\alpha=1$, and so $sR+bR=R$, as required.
\hfill$\Box$

\vskip4mm We have at our disposal all the information necessary to prove:

\vskip4mm \hspace{-1.8em} {\bf Theorem 3.7.}\ \ {\it The following are equivalent for a ring $R$:}
\vspace{-.5mm}
\begin{enumerate}
\item [(1)]{\it $R$ is an elementary divisor ring.}
\vspace{-.5mm}
\item [(2)]{\it $R$ is a B$\acute{e}$zout locally stable ring.}
\vspace{-.5mm}
\item [(3)]{\it $R$ is a B$\acute{e}$zout ring of neat range 1.}
\end{enumerate}
\vspace{-.5mm} {\it Proof.}\ \ $(1)\Rightarrow (3)$ Clearly, $R$ is a B$\acute{e}$zout ring. For any $a,b\in R$ such that $aR+bR=R$, we have $aR+aR+bR=R$.
By virtue of [12, Theorem 32], there exist $y,r,s\in R$ such that $a+by=rs$ and $rR+sR=rR+aR=sR+bR=R$. According to Lemma 3.6, $R/(a+by)R$ is clean. Thus, $R$ is of neat range 1.

$(3)\Rightarrow (2)$ For any $a,b\in R$ such that $aR+bR=R$, there exists a $y\in R$ such that $R/(a+by)R$ is clean. In view of [3, Theorem 17.2.2],
$R/(a+by)R$ has stable range 1, as desired.

$(2)\Rightarrow (1)$ This is obvious by Theorem 3.3.
\hfill$\Box$

\vskip4mm As an immediate consequence of Theorem 3.7, we now derive

\vskip4mm \hspace{-1.8em} {\bf Corollary 3.8.}\ \ {\it Every elementary divisor ring has neat range 1.}

\vskip4mm A ring $R$ is morphic provided that $R/aR\cong r(a)$ for any $a\in R$. For instance, every regular ring is morphic. For any B$\acute{e}$zout domain $R$,
we note that $R/cR$ is morphic for all nonzero $c\in R$ ([16, Theorem 2]).

\vskip4mm \hspace{-1.8em} {\bf Corollary 3.9.}\ \ {\it Let $R$ be a morphic ring. Then the following are equivalent:}
\vspace{-.5mm}
\begin{enumerate}
\item [(1)]{\it $R$ is an elementary divisor ring.}
\vspace{-.5mm}
\item [(2)]{\it $R$ is locally stable.}
\vspace{-.5mm}
\item [(3)]{\it $R$ has neat range 1.}
\end{enumerate}
\vspace{-.5mm} {\it Proof.}\ \ Clearly, $R$ is quasi-morphic. In light of [1, Theorem 15], $R$ is a Bezout ring. Therefore we complete the proof, by
Theorem 3.7.\hfill$\Box$

\vskip4mm In light of Corollary 3.9, every morphic ring of neat range 1 is a ring of stable range 2. Moreover, every morphic elementary divisor ring has
neat range 1. These give affirmative answers of the open problems of Zabavsky and Pihurain [15, Open problems].

\vskip4mm \hspace{-1.8em} {\bf Corollary 3.10.}\ \ {\it Let $R$ be a morphic ring. Then the following are equivalent:}
\vspace{-.5mm}
\begin{enumerate}
\item [(1)]{\it $R$ is an elementary divisor ring.}
\vspace{-.5mm}
\item [(2)]{\it For any $a,b\in R$ such that $r(a)\bigcap r(b)=0$, there exists a $y\in R$ such that $a+by\in R$ is stable.}
\end{enumerate}
\vspace{-.5mm} {\it Proof.}\ \ $(1)\Rightarrow (2)$ For any $a,\in R$ such that $r(a)\bigcap r(b)=0$, we have ${\ell}(r(a)\bigcap r(b))=R$.
In light of [1, Theorem 14], ${\ell}(r(a)\bigcap r(b))=aR+bR$, as $R$ is morphic. Thus, $aR+bR=R$. It follows by Corollary 3.9 that $a+by\in R$ is stable for a $y\in R$, as desired.

$(2)\Rightarrow (1)$ For any $a,b\in R$ such that $aR+bR=R$, then $r(a)\bigcap r(b)=0$. By hypothesis, there exists a $y\in R$ such that $R/(a+by)R$
has stable range 1. This completes the proof, by Corollary 3.9.\hfill$\Box$

\vskip4mm Let $A$ be a ring, and let $E$ be an $A$-module and $R := A\propto E$ be the set of pairs $(a,e)$
with pairwise addition and multiplication given by: $(a,e)(b,f)=(ab, af+be)$.
 We have

\vskip4mm \hspace{-1.8em} {\bf Lemma 3.11.}\ \ {\it Let $A$ be a ring, and let $E$ be an $A$-module. Then $A\propto E$ is locally stable if and only if so is $A$.}
\vskip2mm\hspace{-1.8em} {\it Proof.}\ \ Choose $I=\{
\left(
\begin{array}{cc}
0&b\\
&0
\end{array}
\right)~|~b\in E\}$. Then $I$ is an ideal of $A\propto E$ with $I\subseteq J(A\propto E)$. As $A\propto E/I\cong A$, we obtain the result by Proposition 2.5.\hfill$\Box$

\vskip4mm \hspace{-1.8em} {\bf Theorem 3.12.}\ \ {\it Let $A$ be a domain, $K=qf(A)$, $E$ be a K-vector space, and
$R:=A\propto E$ be the trivial ring extension of $A$ by $E$. Then $R$ is an elementary divisor ring if and only if}
\vspace{-.5mm}
\begin{enumerate}
\item [(1)]{\it $A$ is an elementary divisor ring;}
\vspace{-.5mm}
\item [(2)]{\it $dim_K(E)=1$.}
\end{enumerate}
\vspace{-.5mm} {\it Proof.}\ \ $\Longleftarrow $ Since $R$ is an elementary divisor ring, then $R$ is a B$\acute{e}$zout ring. In view of [2, Theorem 2.1],
$A$ is a B$\acute{e}$zout domain and $dim_K(E)=1$. By virtue of Theorem 3.7, $R$ is locally stable. Hence, $A$ is locally stable, by Lemma 3.11. Thus, $A$ is an elementary divisor ring, by Theorem 3.7.

$\Longrightarrow $ $A$ is an elementary divisor ring, it is a B$\acute{e}$zout domain. In light of [2, Theorem 2.1], $R$ is a B$\acute{e}$zout ring.
In view of Lemma 3.11, $R$ is locally stable. Therefore we complete the proof, by Theorem 3.7.\hfill$\Box$

\vskip4mm As an immediate consequence of Theorem 3.12, we have

\vskip4mm \hspace{-1.8em} {\bf Corollary 3.13.}\ \ Let $A$ be a domain. Then $A\propto qf(A)$ is an elementary divisor ring if and only if so is $A$.

\vskip15mm\bc{\bf 4. Strong Completeness}\ec

\vskip4mm A ring $R$ is strongly completable provided that $a_1R+\cdots +a_nR=dR, a_i,d\in R, i=2,\cdots ,n, $ implies there is a
matrix over $R$ with first row $a_1,\cdots ,a_n$ and $det(A)=d$. One easily checks that a B$\acute{e}$zout ring is strongly completable if and only
if for any $d\in R, d\in
a_1R+\cdots +a_nR, a_i\in R, i=2,\cdots , n, $ implies there is a
matrix over $R$ with first row $a_1,\cdots , a_n$ and $det(A)=d$.

\vskip4mm \hspace{-1.8em} {\bf Theorem 4.1.}\ \ {\it Every locally stable ring is strongly completable.}
\vskip2mm\hspace{-1.8em} {\it Proof.}\ \ Let $R$ be a locally stable ring. Suppose that $a_1R+\cdots +a_nR=dR, a_i,d\in R, i=2,\cdots , n$.
If $n=2$, $d=a_1x_1+a_2x_2$ for some
$x_1,x_2\in R$. Then $-x_2,x_1$ works as a second row. We now consider the assertion for $n\geq 3$. Write $d=a_1x_1+\cdots +a_nx_n, a_1=dq_1,\cdots ,a_n=dq_n$ for some
$x_1,\cdots ,x_n,q_1,\cdots ,q_n\in R$. Let $c=x_1q_1+\cdots +x_nq_n-1$. Then $dc=0$. Further, we have
$q_{1}x_1+\cdots +q_{n-1}x_{n-1}+(q_{n}x_n-c)=1$.
Since $R$ is locally stable, there exists some $x\in R$ such that $S:=R/\big(q_1+q_2x_2x+\cdots +q_{n-1}x_{n-1}x+(q_{n}x_n-c)x\big)R$ has stable range 1.
Clearly,
$\overline{q_{2}}S+\cdots +\overline{q_{(n-1)}}S+\overline{(q_{n}x_n-c)}S=S$. Hence,
we can find a $y_3,\cdots ,y_n\in R$ such that
$$\overline{q_{2}+\cdots +q_3y_3+q_{(n-1)}y_{(n-1)}+(q_{n}x_n-c)y_n}\in U(S).$$ This implies that
$$\big(q_{2}+q_3y_3+q_{(n-1)}y_{(n-1)}+(q_{n}x_n-c)y_n\big)R+\big(q_1+q_2x_2x+\cdots +q_{n-1}x_{n-1}x+(q_{n}x_n-c)x\big)R=R.$$
Clearly, $q_1+q_2x_2x+\cdots +q_{n-1}x_{n-1}x+(q_{n}x_n-c)x
= q_1+\big(q_2+\cdots +q_3y_3+q_{(n-1)}y_{(n-1)}+(q_{n}x_n-c)y_n\big)x_2x+\cdots +q_{n-1}(x_{n-1}x-y_{(n-1)}x_2x+(q_{n}x_n-c)(x-y_nx_2x).$ Thus, we get
$$\big(q_1+q_3s_3+\cdots +q_{n-1}s_{n-1}+(q_{n}x_n-c)s_n\big)R+\big(q_2+q_3t_3+\cdots +q_{n-1}t_{n-1}+(q_{n}x_n-c)t_n\big)R=R.$$
Write $$\big(q_1+q_3s_3+\cdots +q_{n-1}s_{n-1}+(q_{n}x_n-c)s_n\big)s+\big(q_2+q_3t_3+\cdots +q_{n-1}t_{n-1}+(q_{n}x_n-c)t_n\big)t=1.$$
Then $$\left|
\begin{array}{ccccc}
q_1+\sum\limits_{i=3}^{n-1}q_is_i+(q_{n}x_n-c)s_n&q_2+\sum\limits_{i=3}^{n-1}q_it_i+(q_{n}x_n-c)t_n&q_3&\cdots &q_{n}\\
-t&s&&&\\
&&1&&\\
&&&\ddots&\\
&&&&1
\end{array}
\right|\in U(R).$$ By column operations, we see that $$u:=\left|
\begin{array}{ccccc}
q_1-cs_n&q_2-ct_n&q_3&\cdots &q_{n}\\
-t&s&&&\\
**&*&1&&\\
\vdots&\vdots&&\ddots&\\
**&*&&&1
\end{array}
\right|\in U(R).$$
Hence, $$d=u^{-1}d\left|
\begin{array}{ccccc}
q_1-cs_n&q_2-ct_n&q_3&\cdots &q_{n}\\
-t&s&&&\\
**&*&1&&\\
\vdots&\vdots&&\ddots&\\
**&*&&&1
\end{array}
\right|=\left|
\begin{array}{ccccc}
a_1&a_2&a_3&\cdots &a_{n}\\
-tu^{-1}&su^{-1}&&&\\
**&*&1&&\\
\vdots&\vdots&&\ddots&\\
**&*&&&1
\end{array}
\right|,$$ hence the result.\hfill$\Box$

\vskip4mm As immediate consequences of Theorem 4.1, we have

\vskip4mm \hspace{-1.8em} {\bf Corollary 4.2.}\ \ Every ring of neat range 1 is strongly completable.

\vskip4mm \hspace{-1.8em} {\bf Corollary 4.3.} \ \ Every elementary divisor ring is strongly completable.

\vskip4mm Every ring of adequate (neat) range 1 is strongly completable. Thus, every VNL ring is strongly completable.
As every Dedekind domain has almost stable range 1, it follows from Theorem 4.1 that every Dedekind domain is strongly completable [9, Theorem ].

\vskip4mm \hspace{-1.8em} {\bf Example 4.4.}\ \ For any $n\in {\Bbb N}, {\Bbb Z}_n$ is strongly completable. If $n=1$, then $Z_n=Z$ is a Dedekind domain, and so it is strongly completable.
If $n\neq 1$, then $Z_n$ is finite, and so it has stable range 1. Hence, $R$ has almost range 1. In light of Theorem 4.1, $R$ is strongly completable.
This completes the proof.\hfill$\Box$

\vskip4mm An $R$-module $P$ is stably free if $P\oplus R^m\cong R^n$ for some $m,n\in {\Bbb N}$. As is well known, every stably free $R$-module is free if and only
if $R$ is completable, i.e., for any $a_1,\cdots ,a_n\in R$ such that $a_1R+\cdots +a_nR=R$, $(a_1,\cdots ,a_n)$ is the first row of an invertible matrix over
$R$. By the Quillen-Suslin Theorem, $R[x_1,\cdots ,x_n]$ is completable for any principal ideal domain $R$.

\vskip4mm \hspace{-1.8em} {\bf Proposition 4.5.}\ \ {\it Every stably free module over locally stable rings is free.}
\vskip2mm\hspace{-1.8em} {\it Proof.}\ \ Let $R$ be locally stable. Suppose that $a_1R+\cdots +a_nR=R$ for some
$a_1,\cdots ,a_n\in R$. In light of Theorem 4.2, $(a_1,\cdots ,a_n)$ is the first row of an invertible matrix over $R$. That is, $R$ is completable.
Therefore every stably free $R$-module is free.\hfill$\Box$

\vskip4mm \hspace{-1.8em} {\bf Example 4.6.}\ \ Let $R={\Bbb R}[X,Y,Z]$. We claim that $R$ is not locally stable. If not, then ${\Bbb R}[X,Y,Z]/(X^2+Y^2+Z^2-1)$
is locally stable. In view of Proposition 4.5, every stably free $R$-module is free. Let $x,y,z$ denote the images in $R$ of $X,Y,Z$. Note that $(x,y,z)$ is a unimodular column
which defines a split monomorphism $R\to R^3$ with cokernel $P$. Then $P\oplus R\cong R^3$, i.e., $P$ is stably free. In light of [6, Proposition 2.3], $P$ is not free, a contradiction.
Contrast to this assertion, ${\Bbb R}[[X,Y,Z]]$ is locally stable, in terms of Corollary 2.8.

\vskip4mm The preceding example shows that in Proposition 2.7, the condition "$I\subseteq J(R)$" is necessary, as ${\Bbb R}[X,Y,Z]$ is not locally stable while ${\Bbb R}[X,Y,Z]/(X,Y,Z)\cong {\Bbb R}$ is locally stable. In this case $(X,Y,Z)\nsubseteq J\big({\Bbb R}[X,Y,Z]\big)$.

\vskip4mm Let $a\in R$. We use $Z(a)$ denote the set of maximal ideals of $R$ that contains $a$.

\vskip4mm \hspace{-1.8em} {\bf Proposition 4.7.}\ \ {\it Let $R$ be a B$\acute{e}$zout ring. If $Z(a)$ is finite for every $a\not\in J(R)$, then $R$ is strongly completable.}
\vskip2mm\hspace{-1.8em} {\it Proof.}\ \ Suppose that $Z(a)$ is finite for every $a\not\in J(R)$. If $R$ has finitely many distinct maximal ideals, then $R$ is semilocal, and so it has stable range 1. Hence, $R$ is locally stable. If $R$ has infinitely many distinct maximal ideals, then $R$ is an elementary divisor ring, in terms of [5, Corollary 2]. It follows by Theorem 3.7, $R$ is locally stable. Therefore $R$ is locally stable in any case. This completes the proof, by Theorem 4.1.\hfill$\Box$

\vskip15mm \bc{\bf REFERENCES}\ec \vskip4mm {\small
\re{1} V.P. Camillo and W.K. Nicholson, Quasi-morphic rings, {\it J. Algebra Appl.}, {\bf 6}(2007), 789--799.

\re{2} C. Bakkari, On P-B$\acute{e}$zout rings, {\it Internat. J. Algebra}, {\bf 3}(2009), 669--673.

\re{3} H.
Chen, {\it Rings Related Stable Range Conditions}, Series in
Algebra 11, World Scientific, Hackensack, NJ, 2011.

\re{4} O.V. Domsha and I.S. Vasiunyk, Combining local and adequate rings, {\it Book of abstracts of the International Algebraic Conference}, Taras Shevchenko National University of Kyiv, Kyiv, Ukraine, 2014, pp 25.

\re{5} M. Henriksen, Some remarks on elementary divisor rings II,
{\it Michigan Math. J.}, {\bf 3}(1955), 159--163.

\re{6} J.C. McConnell and J.C. Robson, {\it Noncommutative
Noetherian Rings}, New York, NY: John Wiley, 1988.

\re{7} W.W. McGovern, Neat rings, {\it J. Pure Appl. Algebra}, {\bf 205}(2006), 243--265.

\re{8} W.W. McGovern, B\'{e}zout rings with almost stable range
$1$, {\it J. Pure Appl. Algebra}, {\bf 212}(2008), 340--348.

\re{9} M.E. Moore, A strongly complement property of Dedekind
domain, {\it Czechoslovak Math. J.}, {\bf 25}(100)(1975),
282--283.


\re{10} M. Roitman, The Kaplansky condition and rings of almost
stable range $1$, {\it Proc. Amer. Math. Soc.}, {\bf 141}(2013),
3013--3018.

\re{11} B.V. Zabavsky, Diagonal Reduction of Matrices over Rings,
Mathematical Studies Monograph Series, Vol. XVI, VNTL Publisher, 2012.

\re{12} B.V. Zabavsky, Diagonal reduction of matrices over finite stable range rings,
{\it Math. Stud.}, {\bf 41}(2014), 101-108.

\re{13} B.V. Zabavsky and S.I. Bilavska,
Every zero adequate ring is an exchange ring, {\it J. Math. Sci.},
{\bf 187}(2012), 153--156.

\re{14} B.V. Zabavsky and O. Domsha, Diagonalizability theorem for matrices
over certain domains, {\it Algebra $\&$ Discrete Math.}, {\bf 12}(2011), 132--139.

\re{15} B.V. Zabavsky and O.V. Pihura, Conditions under which a morphic ring of stable range 2 is an elementary divisor ring, {\it Book of abstracts of the International Algebraic Conference}, Taras Shevchenko National University of Kyiv, Kyiv, Ukraine, 2014, pp 66.

\re{16} B.V. Zabavsky and O.V. Pihura, B\'{e}zout morphic rings, {\it Visnyk Lviv Univ., Series Mech. Math.}, {\bf 79}(2014), 163--168.
\end{document}